\documentclass[12pt]{article}
\usepackage{mathrsfs}
\usepackage[pdftex]{graphicx}
\usepackage{amsmath,amssymb}
\usepackage{amsfonts}
\usepackage{txfonts}

\newcommand\figcaption{\def\@captype{figure}\caption}

\begin{document}
\title{Diffusive Wave in the Low Mach Limit for Non-Viscous and Heat-Conductive Gas
\thanks{Supported by the Construct Program of the Key Discipline in Hunan Province. email: lyc9009@sina.cn }}
\author{Yechi Liu\\
\small Key Laboratory of High Performance Computing and Stochastic Information Processing\\
\small Department of Mathematics, Hunan Normal University, Changsha, Hunan 410081, P. R. China\\
}
\date{}
\maketitle

\vskip 0.3in
{\bf Abstract}\quad The low Mach number limit for one-dimensional non-isentropic compressible Navier-Stokes system without viscosity is investigated, where the density and temperature have different asymptotic states at far fields. It is proved that the solution of the system converges to a nonlinear diffusion wave globally in time as Mach number goes to zero. It is remarked that the velocity of diffusion wave is proportional with the variation of temperature. Furthermore, it is shown that the solution of compressible Navier-Stokes system also has the same phenomenon when Mach number is suitably small.

\vskip 0.1in
\noindent{\it Keywords:}\quad Compressible Navier-Stokes equations, low Mach limit, nonlinear diffusion wave

\vskip 0.2in

\renewcommand{\theequation}{\thesection.\arabic{equation}}
\section*{1. Introduction}
\setcounter{section}{1}\setcounter{equation}{0}

The one-dimensional compressible Navier-Stokes system in Lagrangian coordinates reads
\begin{eqnarray}
&&v_t-u_x=0,\nonumber\\
&&u_t+P_x=(\mu\frac{u_x}{v})_x,\label{sol1}\\
&&(e+\frac{u^2}{2})_t+(Pu)_x=(\kappa\frac{\theta_x}{v}+\mu\frac{uu_x}{v})_x,\nonumber
\end{eqnarray}
where the unknown functions $v$, $u$ and $\theta$ represent the specific volume, velocity and temperature, respectively, while $\mu>0$ and $\kappa>0$ denote the viscosity and heat conductivity coefficients respectively. Here, we consider the perfect gas, so that the pressure function $P$ and the internal function $e$ are given by
\begin{equation}
P=R\frac{\theta}{v},\quad\textrm{and}\quad e=c_v\theta+const.\nonumber,
\end{equation}
where the parameters $R>0$ and $c_v>0$ are the gas constant and heat capacity at the constant volume respectively. For simplicity, we assume $\mu$ and $\kappa$ are constants, and normalize $R=1$ and $c_v=1$.

The low Mach limit is an important and interesting problem in fluid dynamics. The first result is due to Klainerman and Majda \cite{Kla81,Kla82}, in which they proved the incompressible limit of the isentropic Euler equations to the incompressible Euler equations for local smooth solutions. In \cite{Ala06}, Alazard showed the incompressible limit for Navier-Stokes equations in the whole space. Note that in \cite{Ala06}, the solutions have the same states at the far fields.

Recently, Huang et.al. \cite{Hua16} began to study the case that the solutions have different end states and found that the solutions of compressible Navier-Stokes equations converge to a nonlinear diffusion wave solution globally in time as Mach number goes to zero, which is related to the thermal creep flow. That is, the flow in diffusion wave is only driven by the variation of temperature. This phenomenon is quite different from the constant case. Since the diffusion wave is independent of the viscosity $\mu=0$, we conjecture that the result of \cite{Hua16} is still valid without viscosity, that is $\mu=0$ in the system \eqref{sol1}. Precisely speaking, we consider the non-viscous and heat-conductive gas in the following system
\begin{eqnarray}
&&v_t-u_x=0,\nonumber\\
&&u_t+P_x=0,\label{cha2}\\
&&(\theta+\frac{u^2}{2})_t+(Pu)_x=(\kappa\frac{\theta_x}{v})_x,\nonumber
\end{eqnarray}
where the only difference with the system \eqref{sol1} is that $\mu=0$. We will prove that the solution of system \eqref{cha2} converges to a nonlinear diffusion wave solution globally in time as Mach number tends to zero. Moreover, as the Mach number is suitably small, the flow is only driven by the variation of temperature.

Let $\varepsilon$ be the compressibility parameter, which represents the maximum Mach number of the fluid. As in \cite{Sch07}, we set
\begin{equation}
t\rightarrow\varepsilon t,\quad x\rightarrow x,\quad u\rightarrow\varepsilon u,\quad
\mu\rightarrow\varepsilon\mu,\quad \kappa\rightarrow\varepsilon\kappa.\nonumber
\end{equation}
By the above changes of variables, system \eqref{cha2} is written as
\begin{eqnarray}
&&v_t^\varepsilon-u_x^\varepsilon=0,\nonumber\\
&&u_t^\varepsilon+\frac{1}{\varepsilon^2}P^\varepsilon_x=0,\label{sol3}\\
&&[\theta^\varepsilon+\frac{1}{2}(\varepsilon u^\varepsilon)^2]_t+(P^\varepsilon u^\varepsilon)_x=\kappa(\frac{\theta_x^\varepsilon}{v^\varepsilon})_x.\nonumber
\end{eqnarray}
As $\varepsilon\rightarrow 0$, the limit of solutions of \eqref{sol3} is called the low Mach limit \cite{Kla81,Kla82}. Similar to \cite{Ala06,Hua16,Sch07}, we assume that the pressure is a small perturbation of a given constant state $\bar{P}>0$, i.e.
\begin{equation}\label{PP}
P=\bar{P}+O(\varepsilon),
\end{equation}
and without loss of generality, we further assume $\bar{P}$ to be $1$.
Formally, as the Mach number goes to zero, the limit system of \eqref{sol3} is
\begin{eqnarray}
&&(2u-\kappa T_x)_x=0,\,\,\rho=T^{-1},\nonumber\\
&&u_t+\pi_x=0,\label{bac2}\\
&&\theta_t+u_x=\kappa(\frac{\theta_x}{v})_x,\nonumber
\end{eqnarray}
and
\begin{equation}
v_t-u_x=0,\nonumber
\end{equation}
in which $\pi_x$ is the limit of $\frac{1}{\varepsilon^2}P^\varepsilon_x$.

We will study the low Mach limit and what happens in the limiting process when the background is not constant state, i.e.
\begin{equation}\label{ten}
(v^\varepsilon,\theta^\varepsilon)(x,t)\rightarrow(v_\pm,\theta_\pm),\quad{\rm as}\,x\rightarrow\pm\infty,\quad{\rm with}\,\, \frac{\theta_-}{v_-}=\frac{\theta_+}{v_+},
\end{equation}
where $\theta_-$ may not be equal to $\theta_+$. Following \cite{Hua16}, we shall construct a spacial diffusive wave $(\bar{v},\bar{u},\bar{\theta},\bar{\pi})$ of \eqref{bac2} by choosing
\begin{equation}\label{def1}
(\bar{v},\bar{u},\bar{\theta})\triangleq(T,\frac{\kappa T_x}{2T},T),
\end{equation}
where $T(\eta),\,\eta=\frac{x}{\sqrt{1+t}}$ is the unique self-similar solution of the following diffusion equation
\begin{equation}\label{equ}
T_t=(\frac{\kappa T_x}{2T})_x,\qquad\lim_{x\rightarrow\pm\infty}T=\theta_\pm,
\end{equation}
and $\bar{\pi}$ is a solution of \eqref{bac2}$_2$. For the existence and uniqueness of the self-similar solution of \eqref{equ}, see \cite{Atk74} and \cite{Duy77}. Set $\delta=|\theta_+-\theta_-|$, then the $T(x,t)$ has the following asymptotic expression
\begin{equation}\label{tens}
T_x(x,t)=O(1)\delta(1+t)^{-\frac{1}{2}}e^{-\frac{x^2}{4d(\theta_\pm)(1+t)}},\quad d(\theta_\pm)=\frac{\kappa}{2\theta_\pm},\quad{\rm as}\, x\rightarrow\pm\infty.
\end{equation}
(We will prove this approximation by the idea and method from \cite{Hsi92} in \textbf{Appendix}.)
Because $(\bar{v},\bar{u},\bar{\theta})$ is not the solution of \eqref{sol3}, there will be some non-integrated error terms in the system. It is necessary to introduce a new profile to approximate towards system \eqref{sol3}. That is
\begin{equation}\label{def2}
(\tilde{v},\tilde{u},\tilde{\theta})\triangleq(\bar{v},\bar{u},\bar{\theta}-\frac{1}{2}(\varepsilon\bar{u})^2).
\end{equation}
Then, we have
\begin{equation}
\|(\bar{v}-\tilde{v},\bar{u}-\tilde{u},\bar{\theta}-\tilde{\theta})\|_{L_x^2}\leqslant C\varepsilon^2(1+t)^{-1},\nonumber
\end{equation}
which implies that $(\tilde{v},\tilde{u},\tilde{\theta})$ approximates the diffusive wave solution $(\bar{v},\bar{u},\bar{\theta})$ in $L^2-$norm.
A direct calculation implies
\begin{eqnarray}
&&\tilde{v}_t-\tilde{u}_x=0,\nonumber\\
&&\varepsilon^2\tilde{u}_t+\tilde{P}_x=R_{1x},\label{bac3}\\
&&[\tilde{\theta}+\frac{1}{2}(\varepsilon\tilde{u})^2]_t+(\tilde{P}\tilde{u})_x=\kappa(\frac{\tilde{\theta}_x}{\tilde{v}})_x+R_{2x},\nonumber
\end{eqnarray}
in which
\begin{equation}\label{R1}
R_1=\varepsilon^2\frac{\kappa T_t}{2T}-\varepsilon^2\frac{\tilde{u}^2}{2T}=O(1)\delta\varepsilon^2(1+t)^{-1}e^{-\frac{c_\pm x^2}{1+t}},
\end{equation}
and
\begin{equation}\label{R2}
R_2=\varepsilon^2\kappa\frac{\tilde{u}\tilde{u}_x}{T}-\varepsilon^2\frac{\tilde{u}^3}{2T}
=O(1)\delta\varepsilon^2(1+t)^{-\frac{3}{2}}e^{-\frac{c_\pm x^2}{1+t}},
\end{equation}
as $x\rightarrow\pm\infty$. The initial data of \eqref{sol3} is given by
\begin{equation}\label{indata}
(v^\varepsilon,u^\varepsilon,\theta^\varepsilon)|_{t=0}=(\tilde{v},\tilde{u},\tilde{\theta})(x,0).
\end{equation}
Then we obtain the following global existence and uniform estimates.

\vskip 0.1in

\textbf{Theorem 1 (Uniform Estimates).} Let $(\tilde{v},\tilde{u},\tilde{\theta})$ be the diffusive wave defined by \eqref{def2} and set $\delta=|\theta_+-\theta_-|$ to be the wave strength. Then there exist constants $\varepsilon_0>0$ and $\delta_0>0$, such that if $\varepsilon\leqslant\varepsilon_0$ and $\delta\leqslant\delta_0$, the Cauchy problem \eqref{sol3} with the initial values \eqref{indata} has a unique global smooth solution $(v^\varepsilon,u^\varepsilon,\theta^\varepsilon)$ satisfying
\begin{eqnarray}
\|(v^\varepsilon-\tilde{v},\varepsilon u^\varepsilon-\varepsilon\tilde{u},\theta^\varepsilon-\tilde{\theta})\|_{L_x^2}^2&\leqslant&C\varepsilon^3\sqrt{\delta}(1+t)^{-1+C_0\sqrt{\delta}},\nonumber\\
\|(v^\varepsilon-\tilde{v},\varepsilon u^\varepsilon-\varepsilon\tilde{u},\theta^\varepsilon-\tilde{\theta})_x\|_{L_x^2}^2&\leqslant&C\varepsilon^{\frac{8}{3}}\sqrt{\delta}(1+t)^{-\frac{3}{2}+C_0\sqrt{\delta}},\label{est1}\\
{\rm and}\qquad\quad \|(\theta^\varepsilon-\tilde{\theta})_{xx}\|_{L_x^2}^2&\leqslant&C\varepsilon^{\frac{2}{3}}\sqrt{\delta}(1+t)^{-\frac{3}{2}+C_0\sqrt{\delta}},\nonumber
\end{eqnarray}
where $C$ and $C_0$ are constants independent of $\varepsilon$ and $\delta$, and $\|\cdot\|_{L_x^2}$ denotes the $L^2-$norm with respect to $x$.

\vskip 0.1in

Note that from \eqref{ten}, \eqref{def2}, \eqref{indata} and \eqref{est1}, it is straightforward to see that $(v^\varepsilon-\tilde{v},\varepsilon u^\varepsilon-\varepsilon\tilde{u})\in(L_t^\infty(H_x^1(\mathbb{R})))^2$, $\theta^\varepsilon-\tilde{\theta}\in L_t^\infty(H_x^2(\mathbb{R}))$. Hence the Gagliardo-Nirenberg's inequalitiy
\begin{eqnarray}
\|(v^\varepsilon-\tilde{v},\varepsilon u^\varepsilon-\varepsilon\tilde{u},\theta^\varepsilon-\tilde{\theta})\|_{L_x^{\infty}}^2&\leqslant & C\|(v^\varepsilon-\tilde{v},\varepsilon u^\varepsilon-\varepsilon\tilde{u},\theta^\varepsilon-\tilde{\theta})\|_{L_x^2}\|(v^\varepsilon-\tilde{v},\varepsilon u^\varepsilon-\varepsilon\tilde{u},\theta^\varepsilon-\tilde{\theta})_x\|_{L_x^2}\nonumber\\
{\rm and}\qquad \|(\theta^\varepsilon-\tilde{\theta})_x\|_{L_x^\infty}^2&\leqslant & C\|(\theta^\varepsilon-\tilde{\theta})_x\|_{L_x^2}\|(\theta^\varepsilon-\tilde{\theta})_{xx}\|_{L_x^2},\nonumber
\end{eqnarray}
immediately imply that

\vskip 0.1in

\textbf{Corollary 2 (Low Mach Limit).} Under the assumptions of \textbf{Theorem 1}, in addition, suppose $\delta_0\leqslant\frac{1}{16C_0}$, we have, as $\varepsilon\rightarrow0$,
\begin{eqnarray}
\|(v^\varepsilon-\tilde{v},\theta^\varepsilon-\tilde{\theta})\|_{L_x^\infty}&\leqslant & C\varepsilon^{\frac{17}{12}}\delta^{\frac{1}{4}}(1+t)^{-\frac{1}{2}}\rightarrow0,\nonumber\\
\|u^\varepsilon-\tilde{u}\|_{L_x^\infty}&\leqslant &C\varepsilon^{\frac{5}{12}}\delta^{\frac{1}{4}}(1+t)^{-\frac{1}{2}}\rightarrow0,\label{est2}\\
{\rm and}\quad \|(\theta^\varepsilon-\tilde{\theta})_x\|_{L_x^\infty}&\leqslant &C\varepsilon^{\frac{5}{6}}\delta^{\frac{1}{4}}(1+t)^{-\frac{1}{2}}\rightarrow0.\nonumber
\end{eqnarray}

\vskip 0.1in

From \textbf{Corollary 2}, we shall show the behaviors of solutions of \eqref{sol3} when the Mach number $\varepsilon$ is small.

Without lose of generality, we assume that $\theta_+>\theta_-$. Then for any given constant $\eta_0>0$, by \eqref{tens}, there exists $C_{\eta_0}>0$ such that
\begin{equation}
T^\prime(\eta)>C_{\eta_0}\delta,\quad \textrm{for}\,\, |\eta|\leqslant\eta_0.\nonumber
\end{equation}
From \eqref{def1}, \eqref{def2} and \eqref{est2}, it follows that

\vskip 0.1in

\textbf{Corollary 3 (Driven by the Variation of Temperature).} For the solution obtained by \textbf{Theorem 1}, there exists constants $\varepsilon_1=\varepsilon_1(\eta_0)\leqslant\varepsilon_0$ and $\tilde{C}_0>0$, such that for $\varepsilon\leqslant\varepsilon_1$,
\begin{eqnarray}
&&0<\frac{C_{\eta_0}\delta}{\tilde{C}_0\sqrt{1+t}}<\frac{\tilde{\theta}_x}{\tilde{C}_0}\leqslant u^\varepsilon\leqslant \tilde{C}_0\tilde{\theta}_x,\nonumber\\
&&0<\frac{1}{2}\tilde{\theta}_x\leqslant\theta_x^\varepsilon\leqslant\frac{3}{2}\tilde{\theta}_x,\nonumber
\end{eqnarray}
and
\begin{equation}\label{est4}
\bar{C}_0^{-1}\theta_x^\varepsilon\leqslant u^\varepsilon\leqslant\bar{C}_0\theta_x^\varepsilon, \qquad\textrm{for}\,\,
|x|<\eta_0(1+t)^{\frac{1}{2}}, t\geqslant0,
\end{equation}
where $C_1$ is a positive constant depending only on $\theta_\pm$.

\vskip 0.1in

\textbf{Remark 4.} Due to $\tilde{u}=\frac{\kappa T_x}{2T}$, the velocity $\tilde{u}$ of system \eqref{sol1} is proportional with the variation of temperature $\tilde{\theta}$. Thus the estimate \eqref{est4} shows that, as $\varepsilon$ is suitably small, the velocity $u^\varepsilon$ of system \eqref{sol3} is also proportional with the variation of temperature.

\vskip 0.2in

Finally let us outline the proof of \textbf{Theorem 1}. Since there is no viscosity in \eqref{sol3}$_2$, the system \eqref{sol3} is less dissipative compared with the system \eqref{sol1}. It is not trivial to obtain the higher order estimates. Indeed, even for the basic energy estimate, the derivative estimate for $\psi_y^2$ is missing, see \eqref{ene1} below. To prove \textbf{Theorem 1}, we first obtain the estimates for $\|\zeta_y\|_{H^1}$ though the basic energy. Then, by the compensation method (cf. \cite{Kaw87}), we control simultaneously the estimates for $\|(\phi_y,\psi_y)\|_{H^1}$ by $\|\zeta_y\|_{H^1}$. Motivated by \cite{Hua16}, a new type differential inequality is essentially used to obtain the desired a priori estimates.

The rest of this paper will be arranged in the following way. Section 2 is devoted to the details of the proof, in which the energy estimates will be used, while Section 3 is Appendix to prove an approximation of a self-similar solution, which will be used in the proof of the main results.

\vskip 0.2in

\renewcommand{\theequation}{\thesection.\arabic{equation}}
\section*{2. Proof of Theorem 1}
\setcounter{section}{2}\setcounter{equation}{0}

For simplicity, we omit the superscript $\varepsilon$ of the variables in this section.

\vskip 0.2in
{\bf 2.1 Reformulation of the System}
\vskip 0.1in

Set the scaling
\begin{equation}
y=\frac{x}{\varepsilon},\qquad \tau=\frac{t}{\varepsilon^2},\nonumber
\end{equation}
then, from \eqref{R1} and \eqref{R2}, we have
\begin{equation}\label{R1def}
R_1=O(1)\delta\varepsilon^2(1+\varepsilon^2\tau)^{-1}e^{-\varepsilon^2\frac{c_\pm y^2}{1+\varepsilon^2\tau}},
\end{equation}
and
\begin{equation}\label{R2def}
R_2=O(1)\delta\varepsilon^2(1+\varepsilon^2\tau)^{-\frac{3}{2}}e^{-\varepsilon^2\frac{c_\pm y^2}{1+\varepsilon^2\tau}},
\end{equation}
as $y\rightarrow\pm\infty$.
Define the perturbation around the profile $(\tilde{v},\tilde{u},\tilde{\theta})(\tau,y)$ by
\begin{equation}
(\phi,\psi,\omega,\zeta)(\tau,y)\triangleq(v-\tilde{v},\varepsilon u-\varepsilon \tilde{u},\theta+\frac{1}{2}(\varepsilon u)^2-\tilde{\theta}-\frac{1}{2}(\varepsilon \tilde{u})^2 ,\theta-\tilde{\theta})(\tau,y),\nonumber
\end{equation}
and set
\begin{equation}
(\Phi,\Psi,\bar{W})\triangleq\int_{-\infty}^{y}(\phi,\psi,\omega)(\tau,y)\textrm{d}z.\nonumber
\end{equation}
It is easy to see from \eqref{indata} that $(\Phi,\Psi,\bar{W})(0,\pm \infty)=0$, which enables that $(\Phi,\Psi,\bar{W})$ can be defined in Sobolev Space.

Subtracting \eqref{bac3} from \eqref{sol3} and integrating the resulting system with respect of $y$ yields
\begin{eqnarray}
&&\Phi_\tau-\Psi_y=0,\nonumber\\
&&\Psi_\tau+P-\tilde{P}=-R_1,\label{fin1}\\
&&\bar{W}_\tau+\varepsilon Pu-\varepsilon\tilde{P}\tilde{u}=\kappa(\frac{\theta_y}{v}-\frac{\tilde{\theta}_y}{\tilde{v}})-\varepsilon R_2.\nonumber
\end{eqnarray}
For the variable $\bar{W}$ is the anti-derivative of the total energy, it is more convenient to introduce another variable which is related to the temperature, that is
\begin{equation}
W\triangleq\bar{W}-\varepsilon \tilde{u} \Psi.\nonumber
\end{equation}
By a direct calculation, it holds that
\begin{equation}\label{Ydef}
\zeta=W_y-Y,\qquad{\rm where}\,\,Y\triangleq\frac{1}{2}\Psi_y^2-\varepsilon \tilde{u}_y\Psi.
\end{equation}
Using the new variable $W$, \eqref{fin1} becomes
\begin{eqnarray}
&&\Phi_\tau-\Psi_y=0,\nonumber\\
&&\Psi_\tau-\frac{1}{\tilde{v}}\Phi_y+\frac{1}{\tilde{v}}W_y=Q_1,\label{fin2}\\
&&W_\tau+\Psi_y=\frac{\kappa}{\tilde{v}}W_{yy}+Q_2,\nonumber
\end{eqnarray}
where
\begin{equation}\label{Q1def}
Q_1=J_1+\frac{Y}{\tilde{v}}-R_1,
\end{equation}
\begin{equation}\label{Q2def}
Q_2=\kappa(\frac{1}{v}-\frac{1}{\tilde{v}})\theta_y+J_2-\varepsilon \tilde{u}_\tau\Psi-\frac{\kappa}{\tilde{v}}Y_y-\varepsilon R_2+\varepsilon \tilde{u}R_1
\end{equation}
and
\begin{equation}\label{J1def}
J_1=\frac{\tilde{P}-1}{\tilde{v}}\Phi_y-[P-\tilde{P}+\frac{\tilde{P}}{\tilde{v}}(v-\tilde{v})-\frac{1}{\tilde{v}}(\theta-\tilde{\theta})] =O(1)[|(\Phi_y,W_y)|^2+Y^2+(\varepsilon \tilde{u})^{4}],
\end{equation}
\begin{equation}\label{J2def}
J_2=(1-P)\Psi_y=O(1)[|(\Phi_y,\Psi_y,W_y)|^2+Y^2+(\varepsilon \tilde{u})^4].
\end{equation}
Since the local existence of solution of compressible Navier-Stokes system is known, we only need the following a priori estimates to complete the proof of \textbf{Theorem 1}.

\vskip 0.1in

\textbf{Proposition 5 (A priori Estimates)}. Assume that $(\Phi,\Psi,W)$ is a smooth solution of \eqref{fin2} with zero initial data in the time interval $[0,T]$. Then there exist constants $\varepsilon_1>0$ and $\delta_1>0$ such that if $\varepsilon\leqslant\varepsilon_1$ and $\delta\leqslant\delta_1$, then
\begin{eqnarray}
&&\|(\Phi,\Psi,W)(\tau)\|_{L_y^{\infty}}^2\leqslant C\varepsilon\sqrt{\delta},\nonumber\\
&&\|(\phi,\psi,\zeta)(\tau)\|_{L_y^2}^2\leqslant C\varepsilon^2\sqrt{\delta}(1+\varepsilon^2\tau)^{-1+C_0\sqrt{\delta}},\nonumber\\
&&\|(\phi_y,\psi_y,\zeta_y)(\tau)\|_{L_y^2}^2\leqslant C\varepsilon^{\frac{11}{3}}\sqrt{\delta}(1+\varepsilon^2\tau)^{-\frac{3}{2}+C_0\sqrt{\delta}},\label{est5}\\
&&\|\zeta_{yy}(\tau)\|_{L_y^2}^2\leqslant C\varepsilon^{\frac{11}{3}}\sqrt{\delta}(1+\varepsilon^2\tau)^{-\frac{3}{2}+C_0\sqrt{\delta}}.\nonumber
\end{eqnarray}

\vskip 0.2in

To obtain the above a priori estimates, we need the following a priori assumption.
\begin{equation}\label{ass}
\sup_{0\leqslant\tau\leqslant T} \{\|(\Phi,\Psi,W)(\tau)\|_{L_y^{\infty}}^2+\|(\phi,\psi,\zeta)(\tau)\|_{L_y^2}^2+
\frac{1}{\varepsilon^2}\|(\phi_y,\psi_y,\zeta_y)(\tau)\|_{L_y^2}^2+\frac{1}{\varepsilon^2}\|(\phi_{yy},\psi_{yy},\zeta_{yy})(\tau)\|_{L_y^2}^2\}\leqslant\bar{\delta}^2,
\end{equation}
where $\bar{\delta}$ is a constant depending only on $\delta$.
Immediately, from the a priori assumption and the Galiardo-Nirenberg's inequality, we have
\begin{equation}\label{ass2}
\|(\phi,\psi,\zeta)\|_{L_y^\infty}^2\leqslant C\bar{\delta}.
\end{equation}

\vskip 0.4in

In the sequel, we will use the energy estimates to prove \textbf{Proposition 5}. We will use $\|\cdot\|$ to denote the norm of $L^2$ with respect to $y$ in the rest of this paper.

\vskip 0.2in
{\bf 2.2 Basic Estimates}
\vskip 0.1in

Multiplying \eqref{fin2}$_1$ by $\Phi$, \eqref{fin2}$_2$ by $\tilde{v}\Psi$ and \eqref{fin2}$_3$ by $W$, respectively, and adding all the resultant equations, we have
\begin{equation}\label{ene1}
(\frac{1}{2}\Phi^2+\frac{\tilde{v}}{2}\Psi^2+\frac{1}{2}W^2)_\tau+\frac{\kappa}{\tilde{v}}W_y^2=\ \frac{1}{2}\tilde{v}_\tau\Psi^2-(\frac{\kappa}{\tilde{v}})_yW_yW+\tilde{v}\Psi Q_1+WQ_2+(\Phi\Psi-\Psi W+\frac{\kappa}{\tilde{v}}W_yW)_y.
\end{equation}
Since the order of the term $\tilde{v}\Psi Q_1$ with respect to $(1+t)$ is not enough if we estimate the right-hand side of \eqref{ene1} directly, we need to use a weighted energy method.

Similar to \cite{Hua16}, define
\begin{equation}
m=(\Phi,\Psi,W)^t,
\end{equation}
then \eqref{fin2} can be written as
\begin{equation}\label{jia1}
m_\tau+A_1m_y=A_2m_{yy}+A_3,
\end{equation}
where
\[A_1=\left(\begin{array}{ccc}0&-1&0\\-\frac{1}{\tilde{v}}&0&\frac{1}{\tilde{v}}\\0&1&0\end{array}\right),\quad A_2=\left(\begin{array}{ccc}0&0&0\\0&0&0\\0&0&\frac{\kappa}{\tilde{v}}\end{array}\right),\quad A_3=\left(\begin{array}{c}0\\Q_1\\Q_2\end{array}\right).\]
It is easy to see that $\lambda_1=-\sqrt{\frac{2}{\tilde{v}}}, \lambda_2=0, \lambda_3=-\lambda_1$ are the eigenvalues of $A_1$.

Set
\begin{equation}\label{jia2}
L=\frac{1}{2}(l_1,l_2,l_3)^t,\quad R=\frac{1}{2}(r_1,r_2,r_3),
\end{equation}
where $l_i,r_i,i=1,2,3$ are the left and right corresponding eigenvectors and can be chosen as
\begin{eqnarray}
l_1=(-1,-\frac{2}{\lambda_3},1),&l_2=(\sqrt{2},0,\sqrt{2})&,l_3=(-1,\frac{2}{\lambda_3},1),\nonumber\\
r_1=(-1,-\lambda_3,1)^t,&r_2=(\sqrt{2},0,\sqrt{2})^t&,r_3=(-1,\lambda_3,1)^t.\label{jia0}
\end{eqnarray}
Then, we have
\[l_ir_j=4\delta_{ij},\quad i,j=1,2,3,\qquad\Lambda\triangleq LA_1R=\left(\begin{array}{ccc}\lambda_1&0&0\\0&0&0\\0&0&\lambda_3\end{array}\right),\]
and from \eqref{jia1}, we can get
\begin{equation}\label{jia3}
B_\tau+\Lambda B_y=LA_2RB_{yy}+2LA_2R_yB_y+[(L_\tau+\Lambda L_y)R+LA_2R_{yy}]B+LA_3
\end{equation}
by multiplying \eqref{jia1} by $L$, where
\begin{equation}
B\triangleq Lm\triangleq(b_1,b_2,b_3)^t,\nonumber
\end{equation}
and then
\begin{equation}\label{jia4}
\|(b_1,b_3)\|\sim\|(\Phi,\Psi,W)\|,\quad \|b_2\|\sim\|(\Phi,W)\|.
\end{equation}
As what will be mentioned in \textbf{Appendix}, we could assume $T_y>0$, where $T$ is the self-similar solution of \eqref{equ}, and set
\begin{equation}\label{jia5}
T_1\triangleq\frac{T}{\theta_+},\quad {\rm then}\,\, |T_1-1|\leqslant C\delta.
\end{equation}

Let $N>0$ be a large integer which will be chosen later, and multiply \eqref{jia3} by $\tilde{B}\triangleq(T_1^Nb_1,b_2,T_1^{-N}b_3)$, then we have
\begin{eqnarray}
&&(\frac{T_1^N}{2}b_1^2+\frac{1}{2}b_2^2+\frac{T_1^{-N}}{2}b_3^2)_\tau+
B_yA_4B_y+(\frac{\lambda_1}{2}T_1^Nb_1^2+\frac{\lambda_3}{2}T_1^{-N}b_3^2-\tilde{B}A_4B_y)_y\nonumber\\
&&-\frac{T_1^{N-1}}{2}(N\lambda_1T_{1y}+T_1\lambda_{1y})b_1^2+\frac{T_1^{-N-1}}{2}(N\lambda_3T_{1y}+T_1\lambda_{3y})b_3^2-
\tilde{B}\Lambda L_yRB-(\frac{T_1^N}{2})_\tau b_1^2-(\frac{T_1^{-N}}{2})_\tau b_3^2\nonumber\\
&=&(B-\tilde{B})_yA_4B_y-\tilde{B}A_{4y}B_y+2\tilde{B}LA_2R_yB_y+\tilde{B}(L_\tau R+LA_2R_{yy})B+ \tilde{B}LA_3,\label{jia6}
\end{eqnarray}
where
\[A_4\triangleq LA_2R=\frac{\kappa}{4\tilde{v}}\left(\begin{array}{ccc}1&\sqrt{2}&1\\\sqrt{2}&2&\sqrt{2}\\1&\sqrt{2}&1\end{array}\right)\]
is a nonnegative symmetric matrix.

Set
\begin{eqnarray}
&&E_1\triangleq\int(\frac{1}{2}\Phi^2+\frac{\tilde{v}}{2}\Psi^2+\frac{1}{2}W^2)\textrm{d}y+
\int(\frac{T_1^N}{2}b_1^2+\frac{1}{2}b_2^2+\frac{T_1^{-N}}{2}b_3^2)\textrm{d}y,\nonumber\\
&&\textrm{and}\quad K_1\triangleq\int\frac{\kappa}{\tilde{v}}W_y^2\textrm{d}y+\int B_y A_4 B_y\textrm{d}y,\label{1def}
\end{eqnarray}
then we have $E_1\sim\|(\Phi,\Psi,W)\|^2$ and $\|W_y\|^2\leqslant CK_1$.\\
Since
\begin{equation}
|\tilde{B}\Lambda L_yRB|\leqslant CT_{1y}(T_1^N\lambda_1b_1^2+T_1^{-N}\lambda_3b_3^2),\nonumber
\end{equation}
we can choose $N$ large enough so that
\begin{equation}\label{jia7}
-\frac{T_1^{N-1}}{2}(N\lambda_1T_{1y}+T_1\lambda_{1y})b_1^2+\frac{T_1^{-N-1}}{2}(N\lambda_3T_{1y}-T_1\lambda_{3y})b_3^2-\tilde{B}\Lambda L_yRB\geqslant 2\int|T_y|(b_1^2+b_3^2).
\end{equation}
Next we need to estimate the integral of the right hand side of \eqref{ene1} and \eqref{jia6} term by term.
Using the Cauchy's inequality, we obtain
\begin{equation}
\left|\int\frac{1}{2}\tilde{v}_\tau\Psi^2{\rm d}y\right|\leqslant C\varepsilon^2\delta(1+t)^{-1}\|\Psi\|^2,\nonumber
\end{equation}
and
\begin{equation}
\left|\int(\frac{\kappa}{\tilde{v}})_yW_yW\textrm{d}y\right|\leqslant C\varepsilon^2\delta(1+t)^{-1}\|W\|^2+C\delta\|W_y\|^2.\nonumber
\end{equation}
Note that \eqref{jia4}, then
\begin{equation}\label{jia8}
\left|\int(\tilde{v}\Psi Q_1+WQ_2){\rm d}y\right|\leqslant CI,\quad {\rm where}\,\,I=\int[|(b_1,b_3)|\cdot|(Q_1,Q_2)|+|b_2Q_2|]{\rm d}y.
\end{equation}
For the right-hand side of \eqref{jia6}, from \eqref{jia4}, we have
\begin{equation}
\left|\int[(\frac{T_1^N}{2})_\tau b_1^2+(\frac{T_1^{-N}}{2})_\tau b_3^2]{\rm d}y\right|\leqslant C\varepsilon^2\delta(1+t)^{-1}E_1.\nonumber
\end{equation}
Since $A_4$ is nonnegative and \eqref{jia5} holds, by Cauchy's inequality, we can obtain
\begin{eqnarray}
\left|\int(\tilde{B}-B)_yA_4B_y{\rm d}y\right|&=&\left|\int[(T_1^N-1,0,T_1^{-N}-1)B]_yA_4B_y{\rm d}y\right|\nonumber\\
&\leqslant&\frac{C}{\delta}\int|T_y|^2|B|^2{\rm d}y+C\delta\int|B_y|^2{\rm d}y\nonumber\\
&\leqslant& C\varepsilon^2\delta(1+t)^{-1}E_1+C\delta\|(\Phi_y,\Psi_y,W_y)\|^2.\nonumber
\end{eqnarray}
It is easy to check that the terms $\tilde{B}A_{4y}B_y$, $\tilde{B}LA_2R_yB_y$ and $\tilde{B}(L_\tau R+LA_2R_{yy})B$ satisfy the same estimate.
For the last term, it is obvious from \eqref{jia2} and \eqref{jia0} that
\begin{equation}
\left|\int\tilde{B}LA_3{\rm d}y\right|\leqslant CI,\nonumber
\end{equation}
where $I$ is defined by \eqref{jia8}.\\
Substituting all the above estimates into equations \eqref{ene1} and \eqref{jia6}, using \eqref{1def} and \eqref{jia7}, we have
\begin{equation}
E_{1\tau}+K_1+2\int|T_y|(b_1^2+b_3^2){\rm d}y\leqslant C\varepsilon^2\delta(1+t)^{-1}(E_1+1)+C\delta\|(\Phi_y,\Psi_y,W_y)\|^2+CI.\nonumber
\end{equation}

Next, we only need to calculate $I$ to complete the basic estimate.
Since the term $\varepsilon R_2$ and $\varepsilon\tilde{u}R_1$ in $Q_2$ have better decay rate than the term $R_1$ in $Q_1$, it is naturally that we only need to estimate $\left|\int Q_1b_1{\rm d}y\right|$ and $\left|\int Q_2b_2{\rm d}y\right|$.\\
\eqref{Ydef}, \eqref{J1def} and \eqref{ass2} yield that
\begin{equation}
\int|J_1b_1|{\rm d}y\leqslant C\varepsilon^2\delta(1+t)^{-1}(E_1+1)+C\bar{\delta}\|(\Phi_y,\Psi_y,W_y)\|^2\nonumber
\end{equation}
and
\begin{equation}
\int|\frac{Y}{\tilde{v}}b_1|{\rm d}y\leqslant C\varepsilon^2\delta(1+t)^{-1}E_1+C\bar{\delta}\|\Psi_y\|^2.\nonumber
\end{equation}
From \eqref{R1def}, we have
\begin{equation}
\int|R_1b_1|{\rm d}y\leqslant\delta\int|T_y|b_1^2{\rm d}y+C\varepsilon^2\delta(1+t)^{-1}.\nonumber
\end{equation}
Then, it follows from \eqref{Q1def} that
\begin{equation}
\int|Q_1b_1|{\rm d}y\leqslant\delta\int|T_y|b_1^2{\rm d}y+C\varepsilon^2\delta(1+t)^{-1}(E_1+1)+C\bar{\delta}\|(\Phi_y,\Psi_y,W_y)\|^2.\nonumber
\end{equation}
Similarly, from \eqref{R1def}, \eqref{R2def} and \eqref{J2def}, we obtain
\begin{equation}
\int(|\Phi_y\tilde{\theta}_yb_1+|J_2b_2|+|\varepsilon\tilde{u}_\tau\Phi b_2|+|\varepsilon R_2b_2|+|\varepsilon\tilde{u}R_1b_2|){\rm d}y\leqslant C\varepsilon^2\delta(1+t)^{-1}(E_1+1)+C(\delta+\bar{\delta})\|(\Phi_y,W_y)\|^2.\nonumber
\end{equation}
Using Cauchy's inequality, we obtain from \eqref{ass2} that
\begin{equation}
\int|\Phi_y\zeta_yb_2|{\rm d}y\leqslant C\bar{\delta}(\|\Phi_y\|^2+\|\zeta_y\|^2)\nonumber
\end{equation}
and
\begin{equation}
\int|Y_yb_2|{\rm d}y\leqslant C\varepsilon^2\delta(1+t)^{-1}E_1+C(\delta+\bar{\delta})\|\Psi_y\|^2+C\bar{\delta}\|\psi_y\|^2.\nonumber
\end{equation}
Thus, it follows from \eqref{Q2def} that
\begin{eqnarray}
\left|\int Q_2b_2{\rm d}y\right|&\leqslant&C\int(|\phi\zeta_yb_2|+|\phi\tilde{\theta}_yb_2|+|J_2b_2|+|\varepsilon\tilde{u}_\tau\Phi b_2|+|Y_yb_2|+
|\varepsilon R_2b_2|+|\varepsilon\tilde{u}R_1b_2|){\rm d}y\nonumber\\
&\leqslant&C\varepsilon^2\delta(1+t)^{-1}(E_1+1)+C(\delta+\bar{\delta})\|(\Phi_y,W_y)\|^2+C\bar{\delta}\|(\psi_y,\zeta_y)\|^2.\nonumber
\end{eqnarray}\\
Using the above estimates, we have
\begin{equation}
E_{1\tau}+\frac{1}{2}K_1+\int|T_y|(b_1^2+b_3^2){\rm d}y\leqslant C\varepsilon^2\delta(1+t)^{-1}(E_1+1)+C(\delta+\bar{\delta})\|(\Phi_y,\Psi_y)\|^2+C\bar{\delta}\|(\psi_y,\zeta_y)\|^2\nonumber
\end{equation}
by assuming $\delta$ and $\bar{\delta}$ are suitably small.

Since the norm $\|(\Phi_y,\Psi_y)\|$ cannot be controlled by $K_1$, we need to use the compensation matrix technique.
Multiplying \eqref{fin2}$_2$ by $-\frac{1}{2}\Phi_y$ and \eqref{fin2}$_3$ by $\Psi_y$, respectively, and adding both the resultant equations, we obtain that
\begin{eqnarray}
\frac{1}{2\tilde{v}}\Phi_y^2+\frac{1}{2}\Psi_y^2+(W\Psi_y-\frac{1}{2}\Phi_y\Psi)_\tau\ &=&(W\Psi_\tau-\frac{1}{2}\Phi_\tau\Psi)_y+\frac{1}{2\tilde{v}}W_y\Phi_y-W_y(Q_1+\frac{1}{\tilde{v}}\Phi_y-\frac{1}{\tilde{v}}W_y) \nonumber\\ &&-\frac{1}{2}\Phi_yQ_1+\frac{\kappa}{\tilde{v}}W_{yy}\Psi_y+\Psi_yQ_2,\label{use}
\end{eqnarray}
by using \eqref{fin2}$_1$ and \eqref{fin2}$_2$.
Integrating \eqref{use} with respect of $y$, using \eqref{Ydef} and the Cauchy's inequality, and suppose $\delta$ and $\bar{\delta}$ are suitably small, we can get
\begin{equation}
\int(\frac{1}{4\tilde{v}}\Phi_y^2+\frac{1}{4}\Psi_y^2)\textrm{d}y+\left[\int(W\Psi_y-\frac{1}{2}\Phi_y\Psi)\textrm{d}y\right]_\tau
\leqslant C_1K_1+C_1\|\zeta_y\|^2+C\varepsilon^3\delta(1+t)^{-\frac{3}{2}}.\nonumber
\end{equation}

Choose $\tilde{C}_1$ large enough so that
\begin{eqnarray}
0\leqslant\frac{1}{4}\int(\Phi_y^2+\Psi_y^2){\rm d}y+\frac{1}{2}\tilde{C}_1E_1&\leqslant&
\int(W\Psi_y-\frac{1}{2}\Phi_y\Psi){\rm d}y+\tilde{C}_1E_1\triangleq\tilde{E}_1\nonumber\\
{\rm and}\qquad \frac{1}{4}\tilde{C}_1K_1&\leqslant&\frac{1}{2}\tilde{C}_1K_1-C_1K_1.\nonumber
\end{eqnarray}
Set
\begin{equation}
\tilde{K}_1=\frac{1}{4}\tilde{C}_1K_1+\int(\frac{1}{4\tilde{v}}\Phi_y^2+\frac{1}{4}\Psi_y^2)\textrm{d}y,\qquad
{\rm then}\quad \|(\Phi_y,\Psi_y,W_y)\|^2\sim\tilde{K}_1.\nonumber
\end{equation}
By a direct calculation, we have

\vskip 0.1in

\textbf{Lemma 6 (Basic Estimate)}. If $\delta$ and $\bar{\delta}$ are suitably small, it holds
\begin{equation}\label{stmest1}
\tilde{E}_{1\tau}+\frac{1}{2}\tilde{K}_1+\int|T_y|(b_1^2+b_3^2){\rm d}y\leqslant C\varepsilon^2\delta(1+t)^{-1}(\tilde{E}_1+1)+C\bar{\delta}\|(\psi_y,\zeta_y)\|^2+C_1\|\zeta_y\|^2.
\end{equation}

\vskip 0.2in
{\bf 2.3 First-order Derivative Estimates}
\vskip 0.1in

In order to estimate $\|(\phi,\psi,\zeta)\|$, we need to use the convex entropy.
By applying $\partial y$ to \eqref{fin1}, we have
\begin{eqnarray}
&&\phi_\tau-\psi_y=0,\nonumber\\
&&\psi_\tau+(P-\tilde{P})_y=-R_{1y},\label{fin3}\\
&&\zeta_\tau+\varepsilon Pu_y-\varepsilon\tilde{P}\tilde{u}_y=\kappa(\frac{\theta_y}{v}-\frac{\tilde{\theta}_y}{\tilde{v}})_y+Q_3,\nonumber
\end{eqnarray}
where
\begin{equation}\label{Q3def}
Q_3=\varepsilon\tilde{P}_y\tilde{u}+\frac{1}{2}\varepsilon^2(\tilde{u}^2)_\tau-\varepsilon R_{2y}
\end{equation}
and we used the fact that $\varepsilon u_\tau+P_y=0$ which comes from \eqref{sol3}$_2$.\\
Set
\begin{equation}
F(s)=s-1-\ln s,\nonumber
\end{equation}
then, it is obvious that $F^\prime(1)=0$, and $F(s)$ is strictly convex around $s=1$. Moreover, using Taylor's formula, we obtain that
\begin{equation}\label{Fest}
C_1\phi^2\leqslant F(\frac{v}{\tilde{v}})\leqslant C_2\phi^2,\qquad C_1\zeta^2\leqslant F(\frac{\theta}{\tilde{\theta}}), F(\frac{\tilde{\theta}}{\theta})\leqslant C_2\zeta^2,
\end{equation}
for some positive constants $C_1, C_2>0$.
By a direct calculation, we get
\begin{equation}\label{Fcal1}
[\tilde{\theta}F(\frac{v}{\tilde{v}})]_\tau =\tilde{\theta}_\tau F(\frac{v}{\tilde{v}})-\tilde{\theta}(\frac{1}{v}-\frac{1}{\tilde{v}})\phi_\tau-\frac{\tilde{P}\phi^2}{v\tilde{v}}\tilde{v}_\tau,
\end{equation}
and
\begin{equation}\label{Fcal2}
[\tilde{\theta}F(\frac{\theta}{\tilde{\theta}})]_\tau=\frac{\zeta}{\theta}\zeta_{\tau}-\tilde{\theta}_\tau F(\frac{\tilde{\theta}}{\theta}).
\end{equation}
Multiplying \eqref{fin3}$_2$ by $\psi$ and \eqref{fin3}$_3$ by $\frac{\zeta}{\theta}$, respectively, it holds that
\begin{equation}\label{cal1}
\psi_\tau\psi-\frac{\zeta}{v}\psi_y-\tilde{\theta}(\frac{1}{v}-\frac{1}{\tilde{v}})\phi_y+[(P-\tilde{P})\psi]_y=-\psi R_{1y},
\end{equation}
and
\begin{equation}\label{cal2}
\frac{\zeta}{\theta}\zeta_\tau+\frac{\zeta}{v}\psi_y-\varepsilon\frac{\zeta}{\theta}(P-\tilde{P})\tilde{u}_y =\kappa(\frac{\theta_y}{v}-\frac{\tilde{\theta}_y}{\tilde{v}})_y\frac{\zeta}{\theta}+\frac{\zeta}{\theta}Q_3.
\end{equation}
Therefore, it follows from \eqref{Fcal1}-\eqref{cal2} that
\begin{eqnarray}
[\tilde{\theta}F(\frac{v}{\tilde{v}})+\frac{1}{2}\psi^2+\tilde{\theta}F(\frac{\theta}{\tilde{\theta}})]_\tau&=&\tilde{\theta}_\tau F(\frac{v}{\tilde{v}})-\frac{\tilde{P}\tilde{v}_\tau}{v\tilde{v}}\phi^2-\psi R_{1y}-\tilde{\theta}_\tau F(\frac{\tilde{\theta}}{\theta}) +\varepsilon\frac{\zeta}{\theta}(P-\tilde{P})\tilde{u}_y\nonumber\\&&+ [\kappa(\frac{\theta_y}{v}-\frac{\tilde{\theta}_y}{\tilde{v}})\frac{\zeta}{\theta}-(P-\tilde{P})\psi]_y
-\kappa(\frac{\theta_y}{v}-\frac{\tilde{\theta}_y}{\tilde{v}})(\frac{\zeta}{\theta})_y+\frac{\zeta}{\theta}Q_3.\label{ene2}
\end{eqnarray}

Set
\begin{equation}\label{2def}
E_2=\int[\tilde{\theta}F(\frac{v}{\tilde{v}})+\frac{1}{2}\psi^2+\tilde{\theta}F(\frac{\theta}{\tilde{\theta}})]\textrm{d}y, \qquad \textrm{and}\,\,
K_2=\int\frac{\kappa}{v\theta}\zeta_y^2\textrm{d}y,
\end{equation}
then by \eqref{Fest}, we have $E_2\sim\|(\phi,\psi,\zeta)\|^2$ and there exist a constant $\bar{C}_1>0$ such that
\begin{equation}\label{barC1}
\|\zeta_y\|^2\leqslant\bar{C}_1K_2.
\end{equation}
Next we shall estimate the integral of the terms of the right hand side of \eqref{ene2} as in section 2.2.
It is easy to see from the definition of the function $F$ that
\begin{equation}
\left|\int\tilde{\theta}_\tau F(\frac{v}{\tilde{v}})\textrm{d}y\right|\leqslant C\varepsilon^2\delta(1+t)^{-1}\|\phi\|^2,\nonumber
\end{equation}
and
\begin{equation}
\left|\int\tilde{\theta}_\tau F(\frac{\tilde{\theta}}{\theta})\textrm{d}y\right|\leqslant C\varepsilon^2\delta(1+t)^{-1}\|\zeta\|^2.\nonumber
\end{equation}
Note that $\Psi=\frac{1}{2}\lambda_3(b_3-b_1)$, which together with \eqref{R1def}, yield that
\begin{equation}
\left|\int\psi R_{1y}\textrm{d}y\right|=\left|\int\Psi R_{1yy}\textrm{d}y\right|\leqslant
\varepsilon^2(1+t)^{-1}\int|T_y|(b_1^2+b_3^2){\rm d}y+C\varepsilon^4\delta(1+t)^{-2},\nonumber
\end{equation}
and from \eqref{Q3def}, \eqref{def2} and \eqref{R2def}, we could obtain
\begin{equation}
\left|\int\frac{\zeta}{\theta}Q_3\textrm{d}y\right|\leqslant C\varepsilon^2\delta(1+t)^{-1}\|\zeta\|^2+C\varepsilon^5\delta(1+t)^{-\frac{5}{2}}.\nonumber
\end{equation}
With direct calculation, we have
\begin{equation}
\left|\int\frac{\tilde{P}\tilde{v}_\tau}{v\tilde{v}}\phi^2\textrm{d}y\right|\leqslant C\varepsilon^2\delta(1+t)^{-1}\|\phi\|^2,\nonumber
\end{equation}
and
\begin{equation}
\left|\int\varepsilon\frac{\zeta}{\theta}(P-\tilde{P})\tilde{u}_y\textrm{d}y\right|\leqslant C\varepsilon^2\delta(1+t)^{-1}\|(\phi,\zeta)\|^2.\nonumber
\end{equation}
For the term containing $\kappa$, we need to reform it by
\begin{equation}
\kappa(\frac{\theta_y}{v}-\frac{\tilde{\theta}_y}{\tilde{v}})(\frac{\zeta}{\theta})_y=\frac{\kappa}{v\theta}\zeta_y^2 +\frac{\kappa}{\theta}\tilde{\theta}_y(\frac{1}{v}-\frac{1}{\tilde{v}})\zeta_y +\kappa(\frac{1}{\theta})_y[\frac{\zeta_y}{v}+\tilde{\theta}_y(\frac{1}{v}-\frac{1}{\tilde{v}})]\zeta.\nonumber
\end{equation}
Then, since
\begin{equation}
\left|\int\frac{\kappa}{\theta}\tilde{\theta}_y\zeta_y(\frac{1}{v}-\frac{1}{\tilde{v}})\textrm{d}y\right|\leqslant C\varepsilon^2\delta(1+t)^{-1}\|\phi\|^2+C\delta\|\zeta_y\|^2\nonumber
\end{equation}
and
\begin{equation}
\left|\int\kappa(\frac{1}{\theta})_y\zeta[\frac{\zeta_y}{v}+\tilde{\theta}_y(\frac{1}{v}-\frac{1}{\tilde{v}})]\textrm{d}y\right|\leqslant C(\delta+\bar{\delta}+\sqrt{\bar{\delta}})\|\zeta_y\|^2+C\varepsilon^2\delta(1+t)^{-1}\|(\phi,\zeta)\|^2\nonumber
\end{equation}
from the a priori assumption, the following lemma holds.

\vskip 0.1in

\textbf{Lemma 7 (First-order Derivative Estimate)}. If $\delta$ and $\bar{\delta}$ are suitably small, it holds
\begin{equation}\label{eneest2}
E_{2\tau}+\frac{1}{2}K_2\leqslant\varepsilon^2(1+t)^{-1}\int|T_y|(b_1^2+b_3^2){\rm d}y+C\varepsilon^2\delta(1+t)^{-1}E_2+C\varepsilon^4\delta(1+t)^{-2}.
\end{equation}

\vskip 0.2in

Note that there is the norm of $\psi_y$ in the right hand side of \eqref{stmest1}, we need to control it as what we did to the norm of $\Phi_y,\Psi_y$ in section 2.2. Rewrite \eqref{fin3} in the following form
\begin{eqnarray}
&&\phi_\tau-\psi_y=0,\nonumber\\
&&\psi_\tau-\frac{\theta}{v^2}\phi_y+\frac{1}{v}\zeta_y=J_3,\label{1}\\
&&\zeta_\tau+\frac{\theta}{v}\psi_y=J_4,\nonumber
\end{eqnarray}
where
\begin{equation}
J_3=-(\frac{1}{v}-\frac{1}{\tilde{v}})\tilde{\theta}_y+(\frac{\theta}{v^2}-\frac{\tilde{\theta}}{\tilde{v}^2})\tilde{v}_y-R_{1y},\quad
J_4=-\varepsilon(P-\tilde{P})\tilde{u}_y+\kappa(\frac{\theta_y}{v}-\frac{\tilde{\theta}_y}{\tilde{v}})_y-Q_3.\nonumber
\end{equation}
Multiplying \eqref{1}$_2$ by $-\frac{1}{2}\phi_y$ and \eqref{1}$_3$ by $\frac{v}{\theta}\psi_y$, respectively, and adding both the resultant equations,we can obtain that
\begin{eqnarray}
\frac{\theta}{2v^2}\phi_y^2+\frac{1}{2}\psi_y^2+(\frac{v}{\theta}\psi_y\zeta-\frac{1}{2}\phi_y\psi)_\tau&=&
-\frac{1}{2}(\psi_y\psi)_y+\frac{1}{2}\psi_y^2+\frac{1}{2v}\phi_y\zeta_y+(\frac{v}{\theta})_\tau\psi_y\zeta+\frac{v}{\theta}(\psi_\tau\zeta)_y \nonumber\\ &&-(\frac{1}{v}\phi_y-\frac{1}{\theta}\zeta_y+\frac{v}{\theta}J_3)\zeta_y-\frac{1}{2}\phi_yJ_3+\frac{v}{\theta}\psi_yJ_4.\label{use2}
\end{eqnarray}
Integrating \eqref{use2} with respect of $y$, using the a priori assumption and the Cauchy's inequality, and suppose $\delta$ and $\bar{\delta}$ are suitably small, we can get
\begin{equation}\label{use22}
\int(\frac{\theta}{4v^2}\phi_y^2+\frac{1}{4}\psi_y^2)\textrm{d}y+\left[\int(\frac{v}{\theta}\zeta\psi_y-\frac{1}{2}\phi_y\psi)\textrm{d}y\right]_\tau\leqslant C_2K_2+C_2\|\zeta_{yy}\|^2+C\varepsilon^5\delta(1+t)^{-\frac{5}{2}}.
\end{equation}

\vskip 0.2in
{\bf 2.4 Higher-order Derivative Estimates}
\vskip 0.1in

For the second-order derivative estimate, we need to rewrite \eqref{fin3} in the following form
\begin{eqnarray}
&&\phi_\tau-\psi_y=0,\nonumber\\
&&\psi_\tau-\frac{P}{v}\phi_y+\frac{1}{\tilde{v}}\zeta_y=
(\frac{P}{v}-\frac{\tilde{P}}{\tilde{v}})\tilde{v}_y-(\frac{1}{v}-\frac{1}{\tilde{v}})\theta_y-R_{1y},\label{fin4}\\
&&\zeta_\tau+\psi_y=\kappa(\frac{\theta_y}{v}-\frac{\tilde{\theta}_y}{\tilde{v}})_y+
Q_3+\frac{\varepsilon^2\tilde{u}^2}{2\tilde{v}}\psi_y-\varepsilon(P-\tilde{P})u_y,\nonumber
\end{eqnarray}
according to \eqref{def2}. Then applying $\partial y$ to \eqref{fin4} yields that
\begin{eqnarray}
&&\phi_{y\tau}-\psi_{yy}=0,\nonumber\\
&&\psi_{y\tau}-\frac{1}{\tilde{v}}\phi_{yy}+\frac{1}{\tilde{v}}\zeta_{yy}=Q_4-R_{1yy},\label{fin5}\\
&&\zeta_{y\tau}+\psi_{yy}=\kappa(\frac{\theta_y}{v}-\frac{\tilde{\theta}_y}{\tilde{v}})_{yy}+Q_5,\nonumber
\end{eqnarray}
where
\begin{equation}
Q_4=\frac{\tilde{P}-1}{\tilde{v}}\phi_{yy}+(\frac{P}{v}-\frac{\tilde{P}}{\tilde{v}})\phi_{yy}+
2(\frac{\theta_y v_y}{v^2}-\frac{\tilde{\theta}_y\tilde{v}_y}{\tilde{v}^2})-2(\frac{\theta v_y^2}{v^3}-\frac{\tilde{\theta}\tilde{v}_y^2}{\tilde{v}^3})-
(\frac{1}{v}-\frac{1}{\tilde{v}})\theta_{yy}+(\frac{\theta}{v^2}-\frac{\tilde{\theta}}{\tilde{v}^2})\tilde{v}_{yy},\nonumber
\end{equation}
and
\begin{equation}
Q_5=Q_{3y}+\frac{\varepsilon^2\tilde{u}^2}{2\tilde{v}}\psi_{yy}+\frac{\varepsilon^2}{2}(\frac{\tilde{u}^2}{\tilde{v}})_y\psi_y-\varepsilon(P-\tilde{P})_yu_y -\varepsilon(P-\tilde{P})u_{yy}.\nonumber
\end{equation}
Multiplying \eqref{fin5}$_1$ by $\phi_y$, \eqref{fin5}$_2$ by $\tilde{v}\psi_y$ and \eqref{fin5}$_3$ by $\zeta_y$, respectively, and adding all the resultant equations, we can obtain that
\begin{equation}\label{ene3}
(\frac{1}{2}\phi_y^2+\frac{\tilde{v}}{2}\psi_y^2+\frac{1}{2}\zeta_y^2)_\tau=
(\frac{\tilde{v}}{2})_\tau\psi_y^2+(\phi_y\psi_y-\psi_y\zeta_y)_y+\kappa\zeta_y(\frac{\theta_y}{v}-\frac{\tilde{\theta}_y}{\tilde{v}})_{yy}+\zeta_yQ_5+
\tilde{v}\psi_y(Q_4+R_{1yy}).
\end{equation}

Set
\begin{equation}
E_3=\int(\frac{1}{2}\phi_y^2+\frac{\tilde{v}}{2}\psi_y^2+\frac{1}{2}\zeta_y^2)\textrm{d}y, \qquad \textrm{and}\,\,
K_3=\int\frac{\kappa}{v}\zeta_{yy}^2\textrm{d}y,\nonumber
\end{equation}
then we have $E_3\sim\|(\phi_y,\psi_y,\zeta_y)\|^2$ and there exist a constant $\bar{C}_2>0$ such that
\begin{equation}
\|\zeta_{yy}\|^2\leqslant\bar{C}_2K_3.\nonumber
\end{equation}
Similar to what we did in section 2.2 and 2.3, the integral of the first term in the right hand side of \eqref{ene3} can be estimated by
\begin{equation}
\left|\int(\frac{\tilde{v}}{2})_\tau\psi_y^2\textrm{d}y\right|\leqslant C\varepsilon^2\delta(1+t)^{-1}\|\psi_y\|^2.\nonumber
\end{equation}
Dealing with the term containing $\kappa$, we need to reform it by
\begin{equation}
\kappa\zeta_y(\frac{\theta_y}{v}-\frac{\tilde{\theta}_y}{\tilde{v}})_{yy}=-\frac{\kappa}{v}\zeta_{yy}^2 +\kappa[\zeta_y(\frac{\theta_y}{v}-\frac{\tilde{\theta}_y}{\tilde{v}})_{y}]_y-(\frac{\kappa}{v})_y\zeta_y\zeta_{yy} -\kappa[(\frac{1}{v}-\frac{1}{\tilde{v}})\tilde{\theta}_y]_y\zeta_{yy}.\nonumber
\end{equation}
Then, we have
\begin{eqnarray}
\left|\int(\frac{\kappa}{v})_y\zeta_y\zeta_{yy}\textrm{d}y\right|&\leqslant&
C\left(\int|\phi_y\zeta_y\zeta_{yy}|{\rm d}y+\int|\tilde{v}_y\zeta_y\zeta_{yy}|{\rm d}y\right)\nonumber\\
&\leqslant&(\delta+\sqrt{\bar{\delta}})\|\zeta_{yy}\|^2+C\varepsilon^2\delta(1+t)^{-1}\|\zeta_y\|^2+\|(\phi_y,\zeta_y)\|^{\frac{9}{2}}\nonumber
\end{eqnarray}
by using Gagliardo-Nirenberg's and Young's inequalities and the a priori assumption \eqref{ass}, and
\begin{eqnarray}
\left|\int\kappa[(\frac{1}{v}-\frac{1}{\tilde{v}})\tilde{\theta}_y]_y\zeta_{yy}\textrm{d}y\right|&\leqslant&
\int|\phi_y\tilde{\theta}_y\zeta_{yy}|{\rm d}y+C\int|\phi\tilde{\theta}_{yy}\zeta_{yy}|{\rm d}y\nonumber\\
&\leqslant&\delta\|\zeta_{yy}\|^2+C\varepsilon^2\delta(1+t)^{-1}\|\phi_y\|^2+C\varepsilon^4\delta(1+t)^{-2}\|\phi\|^2.\nonumber
\end{eqnarray}

Note that we only have $\|\zeta_{yy}\|^2$ in the left hand side, it needs to be careful to estimate the rest terms, not to bring the norms of $\phi_{yy}$ and $\psi_{yy}$.\\
Since
\begin{equation}
\left|\int\zeta_yQ_{3y}{\rm d}y\right|\leqslant C\varepsilon\delta(1+t)^{-\frac{1}{2}}\|\zeta_y\|^2+C\varepsilon^8\delta(1+t)^{-4},\nonumber
\end{equation}
\begin{equation}
\left|\int(\frac{\varepsilon^2\tilde{u}^2}{2\tilde{v}}\psi_y)_y\zeta_y{\rm d}y\right|=
\left|\int\frac{\varepsilon^2\tilde{u}^2}{2\tilde{v}}\psi_y\zeta_{yy}{\rm d}y\right|
\leqslant\delta\|\zeta_{yy}\|^2+C\varepsilon^4\delta(1+t)^{-1}\|\psi_y\|^2\nonumber
\end{equation}
and
\begin{eqnarray}
\left|\int[\varepsilon(P-\tilde{P})u_y]_y\zeta_y{\rm d}y\right|&=&\left|\int(P-\tilde{P})(\psi_y+\varepsilon\tilde{u}_y)\zeta_{yy}{\rm d}y\right|\nonumber\\
&\leqslant&\bar{\delta}^{\frac{1}{4}}\|\zeta_{yy}\|^2+C\varepsilon^4\delta(1+t)^{-2}\|(\phi,\zeta)\|^2+C\|(\phi_y,\phi_y,\zeta_y)\|^{\frac{7}{2}},\nonumber
\end{eqnarray}
in which we used the Gagliardo-Nirenberg's and Poincare's inequalities, we have
\begin{eqnarray}
\left|\int\zeta_yQ_5{\rm d}y\right|&\leqslant&(\delta+\bar{\delta}^{\frac{1}{4}})\|\zeta_{yy}\|^2+C\varepsilon\delta(1+t)^{-\frac{1}{2}}\|(\psi_y,\zeta_y)\|^2+
C\varepsilon^4\delta(1+t)^{-2}\|(\phi,\zeta)\|^2\nonumber\\
&&+C\|(\phi_y,\phi_y,\zeta_y)\|^{\frac{7}{2}}+C\varepsilon^8\delta(1+t)^{-4},\nonumber
\end{eqnarray}
while
\begin{equation}
\left|\int\tilde{v}\psi_yR_{1yy}{\rm d}y\right|\leqslant C\varepsilon\delta(1+t)^{-\frac{1}{2}}\|\psi_y\|^2+C\varepsilon^6\delta(1+t)^{-3}\nonumber
\end{equation}
from \eqref{R1def}. Furthermore, it holds
\begin{eqnarray}
\left|\int\tilde{v}\psi_y[(\frac{1}{v}-\frac{1}{\tilde{v}})\theta_{yy}+(\frac{\theta}{v^2}-\frac{\tilde{\theta}}{\tilde{v}^2})\tilde{v}_{yy}]{\rm d}y\right|
&\leqslant&C\int|\psi_y\phi\zeta_{yy}|+|\psi_y\phi\tilde{\theta}_{yy}|+|\psi_y(\phi+\zeta)\tilde{v}_{yy}|{\rm d}y\nonumber\\
&\leqslant&\bar{\delta}^{\frac{1}{4}}\|\zeta_{yy}\|^2+C\|(\phi_y,\psi_y)\|^{\frac{7}{2}}+C\varepsilon\delta(1+t)^{-\frac{1}{2}}E_3+C\varepsilon^3\delta(1+t)^{-\frac{3}{2}}E_2\nonumber
\end{eqnarray}
by using the Poincare's inequality.\\
For the term $\tilde{v}\psi_yQ_4$, by a direct calculation, we obtain
\begin{eqnarray}
\left|\int2\tilde{v}\psi_y(\frac{\theta_yv_y}{v^2}-\frac{\tilde{\theta}_y\tilde{v}_y}{\tilde{v}^2}){\rm d}y\right|&\leqslant&
C\int(|\psi_y\zeta_y\phi_y|+|\phi_y\zeta_y\tilde{v}_y|+|\psi_y\tilde{\theta}_y\phi_y|+|\psi_y\tilde{\theta}_y\tilde{v}_y\phi|){\rm d}y\nonumber\\
&\leqslant&\sqrt{\bar{\delta}}\|\zeta_{yy}\|^2+C\varepsilon\delta(1+t)^{-\frac{1}{2}}E_3
+C\varepsilon^3\delta(1+t)^{-\frac{3}{2}}\|\phi\|^2\nonumber\\
&&+C\|(\phi_y,\psi_y,\zeta_y)\|^{\frac{19}{6}}\nonumber
\end{eqnarray}
and
\begin{equation}
\left|\int2\tilde{v}\psi_y(\frac{\theta v_y^2}{v^3}-\frac{\tilde{\theta}\tilde{v}_y^2}{\tilde{v}^3}){\rm d}y\right|
\leqslant\sqrt{\bar{\delta}}\|\phi_{yy}\|^2+C\varepsilon\delta(1+t)^{-\frac{1}{2}}E_3
+C\varepsilon^3\delta(1+t)^{-\frac{3}{2}}E_2+C\|(\phi_y,\psi_y,\zeta_y)\|^{\frac{7}{2}}\nonumber
\end{equation}
where we assumed that $\bar{\delta}$ suitably small.
On the other hand, from \eqref{def2} and \eqref{fin4}$_1$, it follows that
\begin{equation}
\tilde{v}\psi_y\frac{\tilde{P}-1}{\tilde{v}}\phi_{yy}=\frac{1}{2}\varepsilon^2\tilde{u}^2\phi_{yy}\psi_y=
(\frac{1}{2}\varepsilon^2\tilde{u}^2\phi_y\psi_y)_y-\varepsilon^2\tilde{u}\tilde{u}_y\phi_y\psi_y-\frac{1}{4}\varepsilon^2\tilde{u}^2(\phi_y^2)_\tau.\nonumber
\end{equation}
Similarly,
\begin{equation}
\tilde{v}\psi_y(\frac{P}{v}-\frac{\tilde{P}}{\tilde{v}})\phi_{yy}=(\tilde{v}\psi_y(\frac{P}{v}-\frac{\tilde{P}}{\tilde{v}})\phi_y)_y
-[\tilde{v}(\frac{P}{v}-\frac{\tilde{P}}{\tilde{v}})]_y\phi_y\psi_y-\frac{\tilde{v}}{2}(\frac{P}{v}-\frac{\tilde{P}}{\tilde{v}})(\phi_y^2)_\tau.\nonumber
\end{equation}
Since
\begin{equation}
\left|\int\varepsilon^2\tilde{u}\tilde{u}_y\phi_y\psi_y{\rm d}y\right|+
\left|\int[\tilde{v}(\frac{P}{v}-\frac{\tilde{P}}{\tilde{v}})]_y\phi_y\psi_y{\rm d}y\right|\leqslant C\varepsilon\delta(1+t)^{-\frac{1}{2}}\|(\phi_y,\psi_y)\|^2,\nonumber
\end{equation}
and by assuming $\delta$ and $\bar{\delta}$ suitably small, from \eqref{def2}, \eqref{def1}, \eqref{tens} and \eqref{ass2}, the integral of the terms $-\frac{1}{4}\varepsilon^2\tilde{u}^2(\phi_y^2)_\tau$ and $-\frac{\tilde{v}}{2}(\frac{P}{v}-\frac{\tilde{P}}{\tilde{v}})(\phi_y^2)_\tau$ can be absorbed by $E_{3\tau}$. In fact, Set
\begin{equation}
\bar{E}_3=E_3+\int[\frac{1}{4}\varepsilon^2\tilde{u}^2\phi_y^2+\frac{\tilde{v}}{2}(\frac{P}{v}-\frac{\tilde{P}}{\tilde{v}})\phi_y^2]{\rm d}y,\nonumber
\end{equation}
then, $E_3\sim\bar{E}_3$ from \eqref{ass} as $\bar{\delta}$ suitably small. From \eqref{fin4}, \eqref{PP} and \eqref{def2}, we have
\begin{eqnarray}
&&\left|\int\frac{1}{4}\varepsilon^2(\tilde{u}^2)_\tau\phi_y^2{\rm d}y\right|+
\left|\int[\frac{\tilde{v}}{2}(\frac{P}{v}-\frac{\tilde{P}}{\tilde{v}})]_\tau\phi_y^2{\rm d}y\right|\nonumber\\ &\leqslant&
C\varepsilon^2\delta(1+t)^{-1}E_3+C\varepsilon^3\delta(1+t)^{-\frac{3}{2}}E_2+C\varepsilon^6\delta(1+t)^{-3}+C\|(\phi_y,\psi_y,\zeta_y)\|^{\frac{7}{2}}.\nonumber
\end{eqnarray}
Thus, it holds that
\begin{eqnarray}\label{gai1}
\frac{1}{2}\bar{E}_{3\tau}+\frac{1}{4}K_3&\leqslant&C\varepsilon\delta(1+t)^{-\frac{1}{2}}E_3+C\varepsilon^3\delta(1+t)^{-\frac{3}{2}}E_2+C\varepsilon^6\delta(1+t)^{-3}\nonumber\\
&&+C\|(\phi_y,\psi_y,\zeta_y)\|^{\frac{19}{6}}+\sqrt{\bar{\delta}}\|\phi_{yy}\|^2
\end{eqnarray}
by assuming $\delta$ and $\bar{\delta}$ suitably small.

Since the term $\sqrt{\bar{\delta}}\|\phi_{yy}\|^2$ can not be absorbed by the left hand side, we need more calculation. Multiplying \eqref{fin5}$_2$ by $-\frac{1}{2}\phi_{yy}$ and \eqref{fin5}$_3$ by $\psi_{yy}$, respectively, and adding both the resultant equations, we obtain that
\begin{eqnarray}
\frac{1}{2\tilde{v}}\phi_{yy}^2+\frac{1}{2}\psi_{yy}^2+(\zeta_y\psi_{yy}-\frac{1}{2}\phi_{yy}\psi_y)_\tau &=&(\zeta_y\psi_{y\tau}-\frac{1}{2}\phi_{y\tau}\psi_y)_y-\zeta_{yy}(Q_4-R_{1yy}+\frac{1}{\tilde{v}}\phi_{yy}-\frac{1}{\tilde{v}}\zeta_{yy}) \nonumber\\ &&+\frac{1}{2\tilde{v}}\zeta_{yy}\phi_{yy}-\frac{1}{2}\phi_{yy}(Q_4-R_{1yy})\nonumber\\
&&+[\kappa(\frac{\theta_y}{v}-\frac{\tilde{\theta}_y}{\tilde{v}})_{yy}+Q_5]\psi_{yy},\label{ovo1}
\end{eqnarray}
by using \eqref{fin5}$_1$ and \eqref{fin5}$_2$.
Integrating \eqref{ovo1} with respect of $y$, using the Cauchy's inequality, and suppose $\delta$ and $\bar{\delta}$ are suitably small, we can get
\begin{equation}\label{ovo2}
\int(\frac{1}{4\tilde{v}}\phi_{yy}^2+\frac{1}{4}\psi_{yy}^2)\textrm{d}y+\left[\int(\zeta_y\psi_{yy}-\frac{1}{2}\phi_{yy}\psi_y)\textrm{d}y\right]_\tau
\leqslant C_3K_3+C_3\|\zeta_{yyy}\|^2+C\varepsilon^7\delta(1+t)^{-\frac{7}{2}}.
\end{equation}
Applying $\partial y$ to \eqref{fin5} yields that
\begin{eqnarray}
&&\phi_{yy\tau}-\psi_{yyy}=0,\nonumber\\
&&\psi_{yy\tau}-\frac{1}{\tilde{v}}\phi_{yyy}+\frac{1}{\tilde{v}}\zeta_{yyy}=
(\frac{1}{\tilde{v}})_y\phi_{yy}-(\frac{1}{\tilde{v}})_y\psi_{yy}+Q_{4y}-R_{1yyy},\label{ovo3}\\
&&\zeta_{yy\tau}+\psi_{yyy}=\kappa(\frac{\theta_y}{v}-\frac{\tilde{\theta}_y}{\tilde{v}})_{yyy}+Q_{5y}.\nonumber
\end{eqnarray}
Multiplying \eqref{ovo3}$_1$ by $\phi_{yy}$, \eqref{ovo3}$_2$ by $\tilde{v}\psi_{yy}$ and \eqref{ovo3}$_3$ by $\zeta_{yy}$, respectively, and adding all the resultant equations, we can obtain that
\begin{eqnarray}
(\frac{1}{2}\phi_{yy}^2+\frac{\tilde{v}}{2}\psi_{yy}^2+\frac{1}{2}\zeta_{yy}^2)_\tau&=&(\frac{\tilde{v}}{2})_\tau\psi_{yy}^2+(\phi_{yy}\psi_{yy}-\psi_{yy}\zeta_{yy})_y
-(\frac{\tilde{v}_y}{\tilde{v}})_y\phi_{yy}\psi_{yy}+(\frac{\tilde{v}_y}{\tilde{v}})_y\psi_{yy}^2\nonumber\\
&&+\kappa\zeta_{yy}(\frac{\theta_y}{v}-\frac{\tilde{\theta}_y}{\tilde{v}})_{yyy}+\zeta_{yy}Q_5+\tilde{v}\psi_{yy}(Q_{4y}+R_{1yyy}).\label{ovo4}
\end{eqnarray}

Set
\begin{equation}
E_4=\int(\frac{1}{2}\phi_{yy}^2+\frac{\tilde{v}}{2}\psi_{yy}^2+\frac{1}{2}\zeta_{yy}^2)\textrm{d}y, \qquad \textrm{and}\,\,
K_4=\int\frac{\kappa}{v}\zeta_{yyy}^2\textrm{d}y,\nonumber
\end{equation}
and
\begin{equation}
\bar{E}_4=E_4+\int[\frac{1}{4}\varepsilon^2\tilde{u}^2\phi_{yy}^2+\frac{\tilde{v}}{2}(\frac{P}{v}-\frac{\tilde{P}}{\tilde{v}})\phi_{yy}^2]{\rm d}y.\nonumber
\end{equation}
Then, we have $E_4\sim\|(\phi_{yy},\psi_{yy},\zeta_{yy})\|^2$, $E_4\sim\bar{E}_4$ from \eqref{ass} as $\bar{\delta}$ suitably small, and there exists a constant $\bar{C}_3>0$ such that
\begin{equation}
\|\zeta_{yyy}\|^2\leqslant\bar{C}_3K_4.\nonumber
\end{equation}
Just as what we did in second-order derivative estimate, it holds from \eqref{ovo3} that
\begin{equation}\label{bj1}
\bar{E}_{4\tau}+\frac{1}{2}K_4\leqslant C\varepsilon^2\delta(1+t)^{-1}E_4+C\varepsilon^3\delta(1+t)^{-\frac{3}{2}}E_3+(\delta+\sqrt{\bar{\delta}})E_4+C\varepsilon^8\delta(1+t)^{-4}
\end{equation}
by the aim of the proof of Lemma 4.3 in \cite{Fan15}, which made a good estimate on the term containing $Q_{4y}$.\\
Similar to what we did in section 2.2, choose $\tilde{C}_3$ large enough so that
\begin{eqnarray}
0\leqslant\frac{1}{4}\int(\phi_{yy}^2+\psi_{yy}^2){\rm d}y+\frac{1}{2}\tilde{C}_3\bar{E}_3&\leqslant&
\int(\zeta_y\psi_{yy}-\frac{1}{2}\phi_{yy}\psi_y){\rm d}y+\tilde{C}_3\bar{E}_3\nonumber\\
{\rm and}\qquad \frac{1}{4}\tilde{C}_3K_3&\leqslant&\frac{1}{2}\tilde{C}_3K_3-C_3K_3.\nonumber
\end{eqnarray}
Set
\begin{eqnarray}
\hat{E}_4&=&\frac{1}{2}\tilde{C}_3\bar{E}_3+\frac{1}{2}\int(\zeta_y\psi_{yy}-\frac{1}{2}\phi_{yy}\psi_y){\rm d}y+\bar{E}_4,\nonumber\\
\hat{K}_4&=&\frac{1}{8}\tilde{C}_3K_3+\frac{1}{2}\int(\frac{1}{4\tilde{v}}\phi_{yy}^2+\frac{1}{4}\psi_{yy}^2)\textrm{d}y+\frac{1}{2}K_4,\nonumber
\end{eqnarray}
then we have $\hat{E}_4\sim\bar{E}_3+\bar{E}_4$, $\hat{K}_4\sim\bar{E}_4+K_4$, and
\begin{eqnarray}
\hat{E}_{4\tau}+\frac{1}{2}\hat{K}_4&\leqslant&C\varepsilon\delta(1+t)^{-\frac{1}{2}}\hat{E}_4+C\varepsilon^3\delta(1+t)^{-\frac{3}{2}}E_2+C\varepsilon^6\delta(1+t)^{-3}\nonumber\\
&&+C\|(\phi_y,\psi_y,\zeta_y)\|^{\frac{19}{6}}+C\|(\phi_{yy},\psi_{yy},\zeta_{yy})\|^{\frac{19}{6}},\label{bj2}
\end{eqnarray}
from \eqref{gai1} and \eqref{bj1} by assuming $\delta$ and $\bar{\delta}$ suitably small.

Note that from Young's inequality, we have
\begin{equation}
\varepsilon\delta(1+t)^{-\frac{1}{2}}\hat{E}_4\leqslant C\varepsilon^2\delta(1+t)^{-1}\hat{E}_4+\varepsilon^{\frac{2}{3}}\delta\hat{E}_4,\nonumber
\end{equation}
and by using \eqref{ass}, we obtain
\begin{equation}
\|(\phi_y,\psi_y,\zeta_y)\|^{\frac{19}{6}}+\|(\phi_{yy},\psi_{yy},\zeta_{yy})\|^{\frac{19}{6}}\leqslant \varepsilon\sqrt{\bar{\delta}}\hat{E}_4.\nonumber
\end{equation}
Then, \eqref{bj2} can be written as
\begin{equation}\label{bj0}
\hat{E}_{4\tau}+\frac{1}{2}\hat{K}_4\leqslant C\varepsilon^2\delta(1+t)^{-1}\hat{E}_4+C\varepsilon^3\delta(1+t)^{-\frac{3}{2}}E_2+C\varepsilon^6\delta(1+t)^{-3}+C\varepsilon^{\frac{2}{3}}(\delta+\sqrt{\bar{\delta}})\hat{E}_4.
\end{equation}

Choose $\tilde{C}_2$ large enough so that
\begin{eqnarray}
0\leqslant\frac{1}{4}\int(\phi_y^2+\psi_y^2){\rm d}y+\frac{1}{2}\tilde{C}_2E_2&\leqslant&
\int(\frac{v}{\theta}\zeta\psi_y-\frac{1}{2}\phi_y\psi){\rm d}y+\tilde{C}_2E_2,\nonumber\\
\frac{1}{4}\tilde{C}_2K_2&\leqslant&\frac{1}{2}\tilde{C}_2-C_2K_2,\nonumber\\
{\rm and}\qquad \frac{1}{4}\tilde{C}_2\hat{K}_4&\leqslant&\frac{1}{2}\tilde{C}_2\hat{K}_4-C_2\bar{C}_2K_3-C_3\bar{C}_3K_4.\nonumber
\end{eqnarray}
Set
\begin{eqnarray}
\tilde{E}_4&=&E_2+\int(\frac{v}{\theta}\zeta\psi_y-\frac{1}{2}\phi_y\psi){\rm d}y+\frac{1}{2}\tilde{C}_2\varepsilon^{-\frac{2}{3}}\hat{E}_4\nonumber\\
\tilde{K}_4&=&\frac{1}{4}\tilde{C}_2K_2+\int(\frac{\theta}{4v^2}\phi_y^2+\frac{1}{4}\psi_y^2){\rm d}y+\frac{1}{4}\tilde{C}_2\varepsilon^{-\frac{2}{3}}\hat{K}_4,
\quad {\rm then}\; \hat{E}_4\leqslant C\tilde{K}_4.\label{gai6}
\end{eqnarray}
By a direct calculation, we have

\vskip 0.1in

\textbf{Lemma 8 (Higher-order Derivative Estimate)}. If $\delta$, $\bar{\delta}$ and $\varepsilon$ are suitably small, it holds
\begin{equation}\label{eneest3}
\tilde{E}_{4\tau}+\frac{1}{2}\tilde{K}_4\leqslant
C\varepsilon^2(1+t)^{-1}\int|T_y|(b_1^2+b_3^2){\rm d}y+C\varepsilon^2\delta(1+t)^{-1}\tilde{E}_4+C\varepsilon^4\delta(1+t)^{-2}.
\end{equation}

\vskip 0.2in
{\bf 2.5 Proof of Proposition 5}
\vskip 0.1in

Choose $\tilde{C}_4$ large enough so that
\begin{equation}
\frac{1}{4}\tilde{C}_4\tilde{K}_3\leqslant\frac{1}{2}\tilde{C}_4\tilde{K}_3-C_1\tilde{C}_1K_2.\nonumber
\end{equation}
Set
\begin{equation}\label{gai8}
E_5=\tilde{E}_1+\tilde{C}_3\varepsilon^{-2}\tilde{E}_4,\quad {\rm and}\; K_5=\frac{1}{2}\tilde{K}_1+\frac{1}{2}\tilde{C}_3\varepsilon^{-2}\tilde{K}_4.
\end{equation}
Then, from \eqref{stmest1}, \eqref{barC1} and \eqref{eneest3}, by assuming $\bar{\delta}$ and $\varepsilon$ suitably small, we get
\begin{equation}\label{gai9}
E_{5\tau}+K_5\leqslant C_0\varepsilon^2\sqrt{\delta}(1+\varepsilon^2\tau)^{-1}E_5+C_0\varepsilon^2\delta(1+\varepsilon^2\tau)^{-1}.
\end{equation}
Multiplying \eqref{gai9} by $(1+\varepsilon^2\tau)^{-C_0\sqrt{\delta}}$, integrating the resultant inequality over $(0,\tau)$, we obtain that
\begin{equation}\label{gai10}
(1+\varepsilon^2\tau)^{-C_0\sqrt{\delta}}E_5+\int_0^\tau(1+\varepsilon^2s)^{-C_0\sqrt{\delta}}K_5{\rm d}s\leqslant\sqrt{\delta}.
\end{equation}
Thus
\begin{equation}
E_5\leqslant\sqrt{\delta}(1+\varepsilon^2\tau)^{C_0\sqrt{\delta}}\quad {\rm and}\;
\int_0^\tau K_5{\rm d}s\leqslant\sqrt{\delta}(1+\varepsilon^2\tau)^{C_0\sqrt{\delta}}.\nonumber
\end{equation}
Then
\begin{equation}
\|(\Phi,\Psi,W)\|^2\leqslant\sqrt{\delta}(1+\varepsilon^2\tau)^{C_0\sqrt{\delta}}.\nonumber
\end{equation}

In order to get \eqref{est5}$_2$, we need a better decay rate. Multiplying \eqref{eneest3} by $1+\varepsilon^2\tau$, we have
\begin{eqnarray}
[(1+\varepsilon^2\tau)\tilde{E_4}]_\tau+(1+\varepsilon^2\tau)\tilde{K}_4&\leqslant&
C\varepsilon^2\int|T_y|(b_1^2+b_3^2){\rm d}y+C\varepsilon^2\delta\tilde{E}_4+C\varepsilon^4\delta(1+\varepsilon^2\tau)^{-1}\nonumber\\
&\leqslant&C\varepsilon^2K_5+C\varepsilon^4\delta(1+\varepsilon^2\tau)^{-1},\nonumber
\end{eqnarray}
by using \eqref{jia4}, \eqref{gai6} and \eqref{gai8}. Then,
\begin{equation}\label{gai14}
\tilde{E}_4\leqslant C\varepsilon^2\sqrt{\delta}(1+\varepsilon^2\tau)^{-1+C_0\sqrt{\delta}},\quad {\rm and}\;
\int_0^\tau(1+\varepsilon^2s)\tilde{K}_4{\rm d}s\leqslant C\varepsilon^2\sqrt{\delta}(1+\varepsilon^2\tau)^{C_0\sqrt{\delta}},
\end{equation}
which immediately implies
\begin{eqnarray}
\|(\Phi_y,\Psi_y,W_y,\zeta)\|^2&\leqslant& C\varepsilon^2\sqrt{\delta}(1+\varepsilon^2\tau)^{-1+C_0\sqrt{\delta}}\nonumber\\
\|(\phi_y,\psi_y,\zeta_y)\|^2&\leqslant& C\varepsilon^{\frac{8}{3}}\sqrt{\delta}(1+\varepsilon^2\tau)^{-1+C_0\sqrt{\delta}}\label{gai15}\\
\|(\phi_{yy},\psi_{yy},\zeta_{yy})\|^2&\leqslant& C\varepsilon^{\frac{8}{3}}\sqrt{\delta}(1+\varepsilon^2\tau)^{-1+C_0\sqrt{\delta}}.\nonumber
\end{eqnarray}

Since the time-decay rate of $\|(\phi_y,\psi_y,\zeta_y)\|^2$ and $\|\zeta_{yy}\|^2$ in \eqref{gai15} is less than the one in \eqref{est5}, we need a better estimate.
Multiplying \eqref{gai1} by $(1+\varepsilon^2\tau)^{\frac{3}{2}}$ and integrating the resultant inequality over $(0,\tau)$, we obtain from \eqref{gai14} and \eqref{gai15} that
\begin{eqnarray}
&&(1+\varepsilon^2\tau)^{\frac{3}{2}}\hat{E}_{4\tau}\nonumber\\
&\leqslant&C\int_0^\tau[\varepsilon\delta(1+\varepsilon^2s)\hat{E}_4+\varepsilon^3\delta E_2+
\varepsilon^6\delta(1+\varepsilon^2s)^{-\frac{3}{2}}+(1+\varepsilon^2s)^{\frac{3}{2}}(\|(\phi_y,\psi_y,\zeta_y)\|^{\frac{19}{6}}+\|(\phi_{yy},\psi_{yy},\zeta_{yy})\|^{\frac{19}{6}})]{\rm d}s\nonumber\\
&\leqslant&C\varepsilon^{\frac{11}{3}}\sqrt{\delta}(1+\varepsilon^2\tau)^{C_0\sqrt{\delta}},\nonumber
\end{eqnarray}
where we used
\begin{eqnarray}
&&\int_0^\tau(1+\varepsilon^2s)^{\frac{3}{2}}(\|(\phi_y,\psi_y,\zeta_y)\|^{\frac{19}{6}}+\|(\phi_{yy},\psi_{yy},\zeta_{yy})\|^{\frac{19}{6}}){\rm d}s\nonumber\\
&\leqslant&\varepsilon^2\int_0^\tau(1+\varepsilon^2s)(\|(\phi_y,\psi_y,\zeta_y)\|^2+\|(\phi_{yy},\psi_{yy},\zeta_{yy})\|^2){\rm d}s\leqslant C\varepsilon^4\sqrt{\delta}(1+\varepsilon^2\tau)^{C_0\sqrt{\delta}},\nonumber
\end{eqnarray}
which comes from \eqref{gai15}, \eqref{gai6} and \eqref{gai14}. Thus
\begin{equation}\label{ovo}
\|(\phi_y,\psi_y,\zeta_y)\|^2+\|\zeta_{yy}\|^2\leqslant C\hat{E}_4\leqslant C\varepsilon^{\frac{11}{3}}\sqrt{\delta}(1+\varepsilon^2\tau)^{-\frac{3}{2}+C_0\sqrt{\delta}}.
\end{equation}
By using Gagliardo-Nirenberg's inequality, we obtain
\begin{equation}
\|(\Phi,\Psi,W)\|_{L_y^\infty}^2\leqslant C\|(\Phi,\Psi,W)\|\|(\Phi_y,\Psi_y,W_y)\|\leqslant
C\varepsilon\sqrt{\delta}(1+\varepsilon^2\tau)^{-\frac{1}{4}+\frac{1}{2}C_0\sqrt{\delta}}.\nonumber
\end{equation}
Then the a priori assumption \eqref{ass} is closed and the proof of \textbf{Proposition 5} is completed.

\vskip 0.2in

\renewcommand{\theequation}{\thesection.\arabic{equation}}
\section*{3. Appendix}
\setcounter{section}{3}\setcounter{equation}{0}

In this section, we will prove the approximation of $T_x(x,t)$, where $T$ is the self-similar solution $T(\eta),\,\eta=\frac{x}{\sqrt{1+t}}$ of the equation
\begin{equation}\label{aequ1}
T_t=(\frac{\kappa T_x}{2T})_x,\qquad\lim_{x\rightarrow\pm\infty}T=\theta_\pm,
\end{equation}
when $x$ tends to $\pm\infty$.

\vskip 0.1in

\textbf{Lemma 9} If $T(\eta),\,\eta=\frac{x}{\sqrt{1+t}}$ satisfies \eqref{aequ1}, then
\begin{equation}\label{4}
T_x(x,t)=O(1)\delta(1+t)^{-\frac{1}{2}}e^{-\frac{x^2}{4d(\theta_\pm)(1+t)}},\quad d(T)=\frac{\kappa}{2T},\quad{\rm as}\, x\rightarrow\pm\infty.
\end{equation}
From equation \eqref{aequ1}, by direct calculation, we have
\begin{equation}\label{atpar}
T_t=-\frac{x}{2}(1+t)^{-\frac{3}{2}}T^\prime,
\end{equation}
and
\begin{equation}\label{axpar}
T_x=(1+t)^{-\frac{1}{2}}T^\prime,\qquad T_{xx}=(1+t)^{-1}T^{\prime\prime}.
\end{equation}
Plugging them into \eqref{aequ1}, we obtain
\begin{equation}\label{aequ2}
\frac{\eta}{\kappa}T^\prime+(\frac{T^\prime}{T})^\prime=0.
\end{equation}
Multiplying \eqref{aequ2} by $\frac{T}{T^\prime}$, it follows
\begin{equation}\label{aequ3}
\frac{\eta}{\kappa}T+[ln(\frac{T^\prime}{T})]^\prime=0.
\end{equation}
Then, we have
\begin{equation}\label{a1}
T^\prime(\eta)=\left(\left.\frac{T^\prime}{T}\right|_{\eta=\eta_0}\right)T(\eta)e^{-\int_{\eta_0}^\eta\frac{s}{\kappa}T(s)\textrm{d}s},
\end{equation}
for any given $\eta_0\in \mathbb{R}$. Without loss of generality, we assume $\theta_+>\theta_-$, then $T^\prime(\eta)>0$.

Set $M>0$, then for any $\eta,\eta_0\in(-M,M)$, there exist constants $C_1,C_2,c_1,c_2$ such that
\begin{equation}\label{a2}
C_1e^{-c_1(\eta^2-\eta_0^2)}\leqslant\frac{T^\prime(\eta)}{T^\prime(\eta_0)}\leqslant C_2e^{-c_2(\eta^2-\eta_0^2)}.
\end{equation}
Integrate \eqref{a2} on $(-M,M)$ with respect of $\eta_0$, and suppose $M$ is large enough, then there exist constants $C_3,C_4$ independent of $M$ such that
\begin{equation}\label{app4}
\frac{C_1C_3}{2}\leqslant\frac{T(M)-T(-M)}{T^\prime(\eta_0)}\leqslant C_2C_4.
\end{equation}
Let $M\rightarrow\infty$, we can obtain
\begin{equation}\label{difest}
\frac{1}{C_2C_4}\leqslant\frac{T^\prime(\eta_0)}{\delta}\leqslant\frac{2}{C_1C_3},
\end{equation}
which means $T^\prime(\eta_0)\leqslant C\delta$.\\
On the other hand, we have
\begin{equation}
|T^\prime|\leqslant Ce^{-\frac{\theta_\pm}{2\kappa}\eta^2}
\end{equation}
from \eqref{a1}. Then it follows $T^\prime=O(1)\delta e^{-\frac{\theta_\pm}{2\kappa}\eta^2}$ as $\eta\rightarrow\pm\infty$, which means
\begin{equation}
T_x=O(1)\delta(1+t)^{-\frac{1}{2}}e^{-\frac{\theta_\pm x^2}{2\kappa(1+t)}},\qquad \textrm{as}\,\,x\rightarrow\pm\infty.
\end{equation}

\end{document}